\documentclass{jhrs}
\usepackage{amsmath,epic,eepic,xy,url}
\input xypic
\xyoption{all}
\newtheorem{theorem}{Theorem}[section]
  \newtheorem{proposition}[theorem]{Proposition}
  \newtheorem{lemma}[theorem]{Lemma}

  \theoremstyle{definition}
  \newtheorem{definition}[theorem]{Definition}

  \theoremstyle{remark}
  \newtheorem{remark}[theorem]{Remark}
  
  \numberwithin{equation}{section}

\newcommand{\mb}[1]{{\mathbf #1}}
\newcommand{\mc}[1]{{\mathcal #1}}

\begin{document}
\title{Reparametrizations with given stop data} 
\author{Martin Raussen}
\email{raussen@math.aau.dk}
\homepage{www.math.aau.dk/~raussen}
\address{Department of Mathematical Sciences\\ Aalborg University\\
    Denmark\\Fredrik Bajersvej 7G\\ DK-9220 Aalborg {\O}st}



\begin{abstract}
  In \cite{FR:07}, we performed a systematic investigation of
reparametrizations of continuous paths in a Hausdorff space that
relies crucially on a proper understanding of stop data of a (weakly
increasing) reparametrization of the unit interval. I am grateful to
Marco Grandis (Genova) for pointing out to me that the proof of
Proposition 3.7 in \cite{FR:07} is wrong. Fortunately, the statement of
that Proposition and the results depending on it stay correct. It is
the purpose of this note to provide correct proofs.
\end{abstract}

\received{August 19, 2008}
\revised{September 1, 2008}
\published{June 30, 2009}
\submitted{Ronnie Brown}

\volumeyear{2009}
\volumenumber{4}
\issuenumber{1}
\startpage{1}

\maketitle
\section{Reparametrizations with given stop maps}
To make this note self-contained, we need to include some of the basic
definitions from \cite{FR:07}. The set of all (nondegenerate) closed
subintervals of the unit interval $I=[0,1]$ will be denoted by
$\mc{P}_{[\; ]}(I)=\{[a,b]\mid 0\le a<b\le 1\}$.

\begin{definition}
  \begin{itemize}
  \item A \emph{reparametrization} of the unit interval $I$ is a
    weakly increasing continuous self-map $\varphi :I\to I$ preserving
    the end points.
  \item A \emph{non-degenerate} interval $J\subset I$ is a
    $\varphi$-\emph{stop interval} if there exists a value $t\in I$
    such that $\varphi ^{-1}(t)=J$. The value $t=\varphi(J)\in I$ is
    called a $\varphi$-\emph{stop value}.
  \item The set of all $\varphi$-stop intervals will be denoted as
    $\Delta _{\varphi}\subseteq \mc{P}_{[\; ]}(I)$. Remark that the
    intervals in $\Delta _{\varphi}$ are disjoint and that $\Delta
    _{\varphi}$ carries a natural total order. We let
    $D_{\varphi}:=\bigcup_{J\in \Delta_{\varphi}} J\subset I$ denote
    the \emph{stop set} of ${\varphi}$; and $C_{\varphi}\subset I$ the
    set of all stop values.
  \item The $\varphi$-\emph{stop map} $F_{\varphi}:\Delta _{\varphi}\to
    C_{\varphi}$ corresponding to a reparametrization $\varphi$ is
    given by $F_{\varphi}(J)=\varphi (J)$. 
  \end{itemize}
\end{definition}
It is shown in \cite{FR:07} that $F_{\varphi}$ is an
\emph{order-preserving bijection} between (at most) \emph{countable
  sets}. It is natural to ask (and important for some of the results
in \cite{FR:07}) which order-preserving bijections between such sets
arise from some reparametrization:

To this end, let
\begin{itemize}
\item $\Delta\subseteq \mc{P}_{[\; ]}(I)$ denote a subset of
  \emph{disjoint closed} sub-intervals -- equipped with the natural total
  order;
\item $C\subseteq I$ denote a subset with the same cardinality as
  $\Delta$;
\item $F:\Delta\to C$ denote an order-preserving bijection.
\end{itemize}

I am grateful to the referee for pointing out the following lemma and
its proof:

\begin{lemma}\label{lem:count}
  A subset $\Delta\subseteq \mc{P}_{[\; ]}(I)$ of disjoint closed intervals
 is countable.
\end{lemma}

\begin{proof}
  Given a set $\Delta$ of disjoint non-degenerate closed sub-intervals
  of the unit interval $I$, each will contain rational numbers by
  density. By the axiom of choice, choose for each disjoint
  sub-interval a specific rational number contained in that
  sub-interval. The chosen set $\Delta '\subset \mb{Q}$ of rationals
  is countable as a subset of $\mb{Q}$. Combining an enumeration of
  $\Delta '$ with the bijection between $\Delta '$ and $\Delta$
  mapping each interval to its chosen rational yields an enumeration
  of $\Delta$.
\end{proof}

\begin{proposition}\label{prop:suffnec}
  There exists a reparametrization $\varphi$ with $F_{\varphi}=F$ if
  and only if conditions \emph{(1) - (8)} below are satisfied for intervals
  contained in $\Delta$ and for all $0<z<1$:
  \begin{enumerate}
  \item $\min J=\sup_{J'<J}\max J'\Rightarrow F(J)=\sup_{J'<J}F(J')$;
\item $\max J=\inf_{J<J'}\min J'\Rightarrow F(J)=\inf_{J<J'}F(J')$;
\item $\sup_{J'<z}\max J'=\inf_{z<J''}\min J''\Rightarrow
  \sup_{J'<z}F(J')=\inf_{z<J''}F(J'')$;
\item $\sup_{J'<z}\max J'<\inf_{z<J''}\min J''\Rightarrow
  \sup_{J'<z}F(J')<\inf_{z<J''}F(J'')$;
\item $\inf_{0<J}\min J=0\Rightarrow \inf_{0<J}F(J)=0$;
\item $\inf_{0<J}\min J>0\Rightarrow \inf_{0<J}F(J)>0$;
\item $\sup_{J<1}\max J=1\Rightarrow \sup_{J<1}F(J)=1$;
\item $\sup_{J<1}\max J<1\Rightarrow \sup_{J<1}F(J)<1$.
  \end{enumerate}
\end{proposition}

\begin{proof}
  Conditions (1) -- (3), (5) and (7) are necessary for the stop data of a
  \emph{continuous} reparametrization $\varphi$; (4), (6) and (8) are
  necessary to avoid further stop intervals.

  Given a stop map satisfying conditions (1) -- (8), we construct a
  reparametrization $\varphi _F$ with $F(\varphi_F)=F$ as follows: For
  $t\in D=\bigcup_{J\in\Delta}J$, one has to define: $\varphi(t)=F(J)$
  with $t\in J$. This defines a weakly increasing function $\varphi_F$
  on $D$. Conditions (1) and (2) make sure that this function is
  continuous (on $D$). Condition (3) makes it possible to extend
  $\varphi_F$ uniquely to a weakly increasing continuous function on
  the closure $\bar{D}$: $\varphi_F(z)$ is defined as $\sup_{J'<z}F(J')$ for
  $z=\sup_{J'<z}\max J'$ and/or as $\inf_{z<J''}F(J'')$ for
  $z=\inf_{z<J''}\min J$. By (5) and (7), $\varphi_F(0)=0$ and
  $\varphi_F(1)=1$ if $0,1\in \bar{D}$; if not, we have to take these as a
  definition.

  The complement $O=I\setminus \bar{D}$ is an open (possibly empty)
  subspace of $I$, hence a union of at most countably many open
  subintervals $J=[a^J_-,a^J_+]$ with boundary in $\partial D\cup \{
  0,1\}$. Condition (4), (6) and (8) make sure, that
  $\varphi_F(a^J_-)<\varphi_F(a^J_+)$. Hence, every collection of
  strictly increasing homeomorphisms between $[a_-^J,a_+^J]$ and
  $[\varphi_F(a_-^J),\varphi_F(a_+^J)]$ -- preserving endpoints -- extends
  $\varphi_F$ to a continuous increasing map $\varphi_F:I\to I$ with
  $\Delta_{\varphi_F}= \Delta, C_{\varphi_F}=C$ and $F_{\varphi_F}=F$.
\end{proof}

It is natural to ask, whether
\begin{itemize}
\item every at most countable subset $C\subset I$ occurs as set of
  stop values of some reparametrization: This is answered
  affirmatively in \cite{FR:07}, Lemma 2.10;
\item every set $\{ I\}\neq\Delta\subset\mc{P}_{[\; ]}(I)$ of closed
  disjoint intervals arises as set of stop intervals of a
  reparametrization:
\end{itemize}

\begin{proposition}\label{prop:repDelta} 
  For every $\{ I\}\neq\Delta$ of closed disjoint
  sub-intervals in the unit interval $I$, there exists a reparametrization
  $\varphi$ with $\Delta_{\varphi}=\Delta$.
\end{proposition}

\begin{proof} 
  We use Lemma \ref{lem:count} to provide us with an enumeration $j$
  of the totally ordered set $\Delta$ (defined either on $\mathbf{N}$
  or on a finite integer interval $[1,n]$). Using $j$, we are going to
  construct a reparametrization $\varphi$ with stop value set
  $C_{\varphi}$ included in the set $I[\frac{1}{2}]=\{0\le
  \frac{l}{2^k}\le 1\}$ of rational numbers with denominators a power
  of 2. To this end, we will associate to every number $z\in
  I[\frac{1}{2}]$ either an interval in $\Delta$ or a degenerate one
  point interval; we end up with an ordered bijection beween
  $I[\frac{1}{2}]$ and a superset of $\Delta$; all excess intervals
  will be degenerate one-point sets.

  To get started, let $I_0$ denote either \emph{the} interval in
  $\Delta$ containing $0$ or, if no such interval exists, the
  degenerate interval $[0,0]=\{ 0\}$; likewise define $I_1$. Every
  number $z\in I[\frac{1}{2}]$ apart from $0$ and $1$ has a unique
  representation $z=\frac{l}{2^k}$ with $l$ \emph{odd}, $0<l<2^k$. The
  construction proceeds by induction on $k$ using the enumeration $j$.

  Assume for a given $k\ge 1$, $I_z$ and thus the map $I:z\mapsto I_z$
  defined for all $z=\frac{l}{2^{k-1}},\; 0\le l\le 2^{k-1}$ as an
  ordered map. For $0<z=\frac{l}{2^k}<1$ and $l$ odd, both $z_{\pm}
  =z\pm \frac{1}{2^k}$ have a representation as fraction with
  denominator $2^{k-1}$ and thus $I_{z_-}<I_{z_+}$ are already
  defined. Let $I_z=j(m)$ with $m$ minimal such that
  $I_{z_-}<j(m)<I_{z_+}$ if such an $m$ exists; if not, then $I_z$ is
  defined as the degenerate interval containing the single element
  $\frac{1}{2}(\max I_{z_-}+\min I_{z_+})$. The map $I:z\mapsto I_z$
  thus constructed on $I[\frac{1}{2}]$ is order-preserving. Moreover,
  this map is onto, since -- by an induction over $k$ -- $I_{j(k)}$
  occurs as $I_z$ with some $z$ of the form $\frac{l}{2^k}$. Hence,
  there is an order-preserving inverse map $I^{-1}:I_z\mapsto z$.

  For $k\ge 0$, let $\varphi_k$ denote the piecewise linear
  reparametrization that has constant value $z$ on $I_z$ for
  $z=\frac{l}{2^k},\; 0\le l\le 2^k$ and that is linear inbetween
  these intervals. Remark that $\varphi_{k+1}=\varphi_k$ on all $I_z$
  with $z=\frac{l}{2^k}$ including all occuring degenerate
  intervals. As a consequence, $\parallel
  \varphi_k-\varphi_{k+1}\parallel <\frac{1}{2^k}$, and hence for all
  $l>k$, $\parallel \varphi_k-\varphi_l\parallel
  <\frac{1}{2^{k-1}}$. Hence, the sequence
  $(\varphi_k)_{k\in\mathbf{N}}$ converges uniformly to a continuous
  reparametrization $\varphi$.

  By construction, the resulting reparametrization $\varphi$ is
  constant on all intervals in $\Delta$; on every open interval between
  these stop intervals, it is linear and strictly increasing. In
  particular, $\Delta_{\varphi}=\Delta$.
\end{proof}
\begin{remark}
  I was first tempted to prove Proposition \ref{prop:repDelta} by
  taking some integral of the characteristic function of the
  complement of $D$ and to normalize the resulting function. But in
  general, this does not work out since, as already remarked in
  \cite{FR:07}, it may well be that $\bar{D}=I$!
\end{remark}

\section{Concluding remarks}
\begin{remark}
  \begin{enumerate}
  \item Instead of constructing the reparametrization $\varphi$ in
    Proposition \ref{prop:repDelta}, it is also possible to apply the
    criteria in Proposition \ref{prop:suffnec} to the restriction
    $I_{|_{\Delta}}$ of the map $I$ from the proof above.
  \item Proposition \ref{prop:suffnec} replaces Proposition 2.13 in
    \cite{FR:07}. To get sufficiency, requirements (1) and (2) had to
    be added to those mentioned in \cite{FR:07} in order to make sure
    that the map $\varphi_F$ is continuous on $D$. Moreover, (6) and (8) had
    to be added to avoid stop intervals containing $0$, resp.~$1$ in
    case $\Delta$ does not contain such intervals. 

    In particular, the midpoint map $m$ that associates to every
    interval in $\Delta$ its midpoint satisfies the criteria given in
    \cite{FR:07}, Proposition 2.13, but if fails in general to satisfy
    conditions (1) and (2) in Proposition \ref{prop:suffnec} in this
    note; in particular, the map $\varphi_m$ will in general not be
    continuous, as remarked by M.~Grandis. The midpoint map $m$ was
    used in the flawed proof of \cite{FR:07}, Proposition 3.7 --
    stated as Proposition \ref{prop:reg} below.
  \end{enumerate}
\end{remark}

The main focus in \cite{FR:07} is on reparametrizations of continuous
paths $p:I\to X$ into a Hausdorff space $X$. A continuous path $q$ is
called \emph{regular} if it is constant or if the restriction $q{|_J}$
to every non-degenerate sub-interval $J\subseteq I$ is
\emph{non-constant}.

\begin{proposition}\label{prop:reg} \emph{(Proposition 3.7 in} 
\emph{\cite{FR:07})}\\
  For every path $p:I\to X$, there exists a regular path $q$ and a
  reparametrization such that $p=q\circ \varphi$.
\end{proposition}

\begin{proof}
  A non-constant path $p$ gives rise to the set of all (closed
  disjoint) \emph{stop intervals} $\Delta_p\subset \mc{P}_{[\; ]}(I)$,
  consisting of the maximal subintervals $J\subset I$ on which $p$ is
  constant. Proposition \ref{prop:repDelta} yields a reparametrization
  $\varphi$ with $\Delta_{\varphi}=\Delta_p$ and thus a set-theoretic
  factorization
\begin{equation*}
  \xymatrix{%
    I \ar[r]^p \ar[d]_{\varphi} & X \\
    I \ar@{-->}[ur]_{q}
  }\end{equation*}
through a map $q:I\to X$ that is not constant on any non-degenerate 
subinterval $J\subseteq I$. The continuity of $q$ follows 
as in the remaining lines of the proof in \cite{FR:07}.
\end{proof}

\noindent See also the references in \cite{FR:07}.
\end{document}